\documentclass[11pt]{article}
\usepackage{sectsty,latexsym,amssymb,dsfont,textcomp}
\usepackage[centertags]{amsmath}
\usepackage{fullpage}
\usepackage{color}
\usepackage{epsfig,psfrag}

\usepackage[centertags]{amsmath}
\sectionfont{\normalsize}
\subsectionfont{\normalsize}
\begin{document}
\date{}
\numberwithin{equation}{section}
\title{The problem of dynamic cavitation in nonlinear elasticity\\\vspace{6pt}
{ \footnotesize (In: {\it S\'eminaire Laurent Schwartz - EDP et applications (2012-2013)}, Exp. 14, 1-17. DOI: 10.5802/slsedp.41.)\vspace{5pt}}
}

\author{
Jan Giesselmann\footnote{Weierstrass Institute, Berlin, Germany,
     jan.giesselmann@wias-berlin.de}
\and
Alexey Miroshnikov\footnote{Department of Mathematics and Statistics, University of
Massachusetts Amherst, USA,  amiroshn@gmail.com}
\and
Athanasios E. Tzavaras\footnote{Department of Applied Mathematics, University of
Crete, Heraklion, Greece and Institute for Applied and Computational Mathematics, FORTH\@, Heraklion, Greece,
tzavaras@tem.uoc.gr}
}

\maketitle

\newtheorem{lemma}{Lemma}[section]
\newtheorem{theorem}[lemma]{Theorem}
\newtheorem{proposition}[lemma]{Proposition}
\newtheorem{maintheorem}[lemma]{Main Theorem}
\newtheorem{corollary}[lemma]{Corollary}
\newtheorem{definition}[lemma]{Definition}
\newtheorem{remark}[lemma]{Remark}
\newtheorem{remarks}[lemma]{Remarks}
\newtheorem{Notation}[lemma]{Notation}
\newcommand{\proof}{\noindent {\it Proof}\;\;\;}
\newcommand{\qed}{\protect~\protect\hfill $\Box$}

\newcommand{\wkarr}{\; \rightharpoonup \;}
\def\Weak{\,\,\relbar\joinrel\rightharpoonup\,\,}

\font\msym=msbm10
\def\charf {\mbox{{\text 1}\kern-.24em {\text l}}}

\def\A{\mathbb A}
\def\Z{\mathbb Z}
\def\K{\mathbb K}
\def\J{\mathbb J}
\def\L{\mathbb L}
\def\D{\mathbb D}
\def\cD{\mathcal D}
\def\cO{\mathcal O}
\def\cQ{\mathcal Q}
\def\Mink{{\mathop{\hbox{\msym \char '115}}}}
\def\Integers{{\mathop{\hbox{\msym \char '132}}}}
\def\Complex{{\mathop{\hbox{\msym\char'103}}}}
\def\C{\Complex}
\font\smallmsym=msbm7

\newcommand{\del}{\partial}
\newcommand{\eps}{\varepsilon}
\newcommand{\cof}{\hbox{cof}\,}
\newcommand{\dt}{\hbox{det}\,}

\newcommand{\cP}{\mathcal{P}}
\newcommand{\cS}{\mathcal{S}}
\newcommand{\cM}{\mathcal{M}}
\renewcommand{\div}{\hbox{div}\,}

\newcommand{\RR}{\mathbb{R}}
\def\Real{{\mathbb{R}}}
\def\R{\Real}
\def\torus{{\mathbb{T}}}
\def\T{\torus}
\def\charf {{{\text{\rm 1}}\kern-.24em {\text{\rm l}}}}

\newcommand{\BBR}[1]{\left( #1 \right)}

\def\div{\hbox{div}\,}
\def\supp{\hbox{supp}\,}
\def\dist{\hbox{dist}\,}
\newcommand{\tcb}{\textcolor{blue}}
\newcommand{\tcg}{\textcolor{green}}
\newcommand{\tcr}{\textcolor{red}}

%
%

%

\begin{abstract}
\noindent
\noindent
The notion of singular limiting induced from continuum solutions (slic-solutions) is  applied to
the problem of cavitation in nonlinear elasticity, in order to re-assess an example of non-uniqueness of
entropic weak solutions  (with polyconvex energy) due to a forming cavity.
\end{abstract}

%
%
%
\section{Introduction}
\label{intro}
\setcounter{equation}{0}

The equations describing radial motions of isotropic elastic materials admit a special weak solution
describing  a cavity that emerges at a point from a homogeneously deformed state.
Cavitating solutions are self-similar in nature and were introduced by Ball \cite{ball82} and
by  Pericak-Spector and Spector \cite{ps88,ps97}.
They turn out to decrease the total mechanical energy and provide a striking example of non-uniqueness
of entropy weak solutions for polyconvex energies \cite{ps88}.

There is a class of problems in solid mechanics, such as fracture, cavitation or shear bands,
where discontinuous motions emerge from smooth motions via a mechanism of material instability.
Any attempt to study such solutions that lie at the limits of continuum modeling
needs to reckon with the problem of giving a proper definition for such solutions.
Once the material breaks or a shear band forms the motion can no longer be described at the level
of continuum modeling and microscopic modeling or higher-order regularizing mechanisms have to be taken into account.
Still, as such structures develop there is expected an intermediate time scale where both types of modeling apply.
In these lecture notes we present  the notion of singular limiting induced from continuum solution
(or slic-solution)  \cite{GT13}, stating that a discontinuous motion is a slic-solution  if its averages
are a family of smooth approximate solutions to the problem.

The lecture notes are based on two recent works on cavity formation \cite{MT13} and \cite{GT13}.
\cite{MT13} complements \cite{ps88,ps97} and establishes various further properties of weak solutions describing
cavitation, in particular indicating that cavity formation is necessarily associated with a unique precursor shock.
In \cite{GT13} the notion of singular limiting induced from continuum solution is introduced and applied to the
problem of cavity formation. It is shown that there  is an energetic cost for creating the cavity,
which is captured by slic-solutions but neglected by the usual entropic weak solutions.
Once this surface energy cost
is accounted for, the paradox of nonuniqueness is removed, in the sense that the cavitating solution together with
the surface energy of the cavity has higher energy than the homogeneous deformation.

The structure of these  notes is as follows: In section \ref{sec-eqel} we present the system of elasticity and the
general requirements imposed  by mechanical considerations and discuss various notions of relaxed convexity
hypotheses associated to elasticity. In section \ref{sec-radel} we introduce the equations of radial elasticity applicable to
isotropic elastic materials. In section \ref{sec-cavi} we present the problem of cavitation and  outline the construction
and properties of weak solutions describing cavity formation from a homogeneously deformed state.
In section \ref{sec-frac} we introduce an illuminating example for fracture in 1-d and in section \ref{sec-slic} the concept of singular limiting
induced from continuum solutions which gives an interpretation to this example as a solution of the equations of
one-dimensional elasticity. It is shown that there is an additional contribution from the opening crack to the resulting
energy of the fracturing solution. In section \ref{sec-sliccav}, the notion of slic-solution is
extended for the cavitation problem and the effect of this concept to the energy of the cavitating solution is calculated.
It turms out that the slic-solution provides  a more discriminating concept of solution for strong singularities.


\section{The equations of elasticity}
\label{sec-eqel}

The equations of elasticity is the system of nonlinear partial differential equations
\begin{equation}
\label{elas-m}
\frac{\del^2 y }{\del t^2} = \div \frac{\del W}{\del F} (\nabla y)
\end{equation}
where $y(x,t) : \R^d \times \R_+ \to \R^d$ describes the motion. It is customary to introduce the
velocity $v = \frac{\del y}{\del t}$ and the deformation gradient  $F = \nabla y$ and to write \eqref{elas-m}
as a first order hyperbolic system
\begin{equation}
\label{elassys}
\left \{ \quad
\begin{aligned}
\del_t F_{i \alpha} &= \del_\alpha v_i
\\
\del_t v_i &= \del_\alpha \frac{\del W}{\del F_{i \alpha}} (F)
\\
\del_\alpha F_{i \beta} &- \del_\beta F_{i \alpha} = 0
\end{aligned}
\right .
\end{equation}
The last equation stands for a constraint that ensures that $F$ is a gradient. This constraint
does not cause problems as it is an involution \cite{Dafermos86}, that is if it is satisfied for the initial data
then it is propagated by  \eqref{elassys}$_1$ to the solutions. In  \eqref{elas-m} the
hypothesis of hyperelasticity, $S = \frac{\del W}{\del F} $, is used stating that the Piola-Kirchhoff stress $S$
 is derived as a gradient of a stored energy function $W : M^{d\times d} \to \R$.
 This hypothesis makes the theory  consistent with the second law of thermodynamics.

\begin{figure}
  \begin{center}
    \includegraphics[height=4cm,width=8cm]{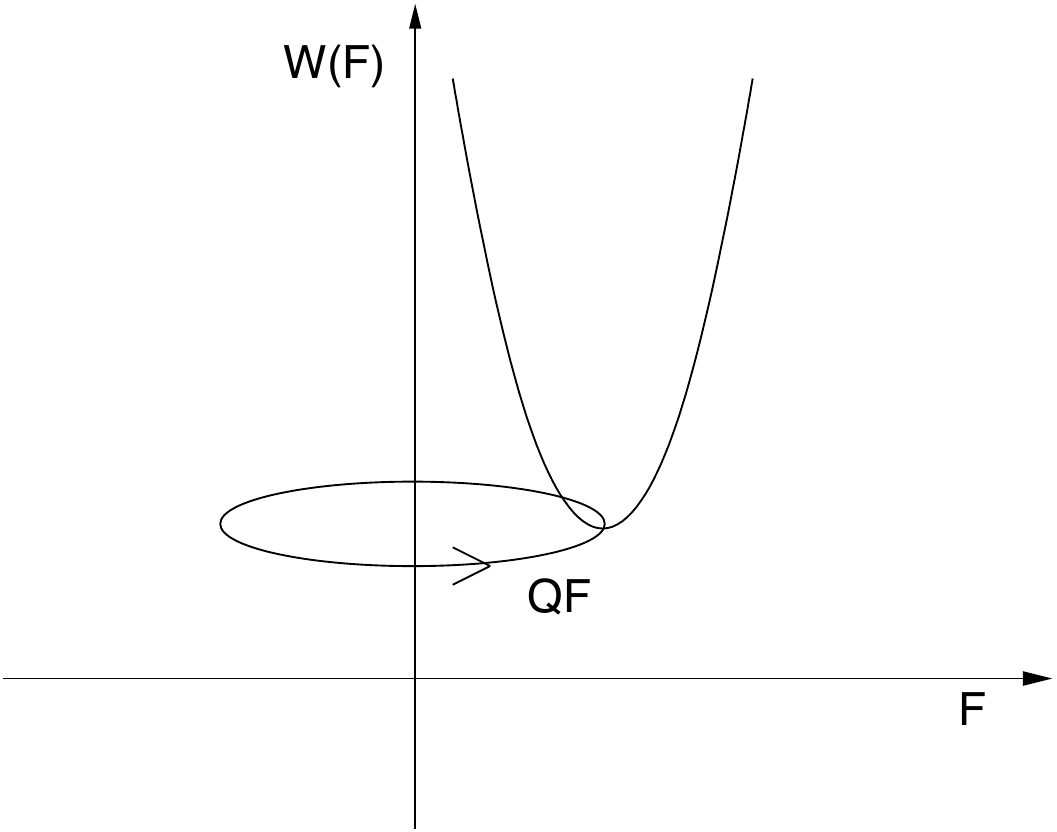}
  \end{center}
  \caption[subfigure]{A depiction of the stored energy}
   \label{figure1}
\end{figure}

The principle of frame indifference dictates that the stored energy  satisfies the invariance
$$
 W (Q F) = W(F) \quad \forall \; Q \in SO(d)
$$
for any matrix $Q$ describing a proper rotation.
In order to interpret a geometric map $y$ as a physical motion, $y$ has to be globally one-to-one
so that matter does not interpenetrate. A necessary condition for achieving that is to cost infinite energy
whenever
a finite volume is compressed down to a zero volume, what dictates the natural assumption for the
stored energy
$$
W(F) \to \infty \quad \mbox{as} \quad  \dt F \to 0 \, .
$$
These two requirements are in general inconsistent with convexity of the stored energy as can be
seen in the schematic picture of Figure \ref{figure1}.
It is thus too restrictive to take $W$ as a convex function.

The system \eqref{elassys} is a system of conservation laws. It is equipped for smooth solutions with the
 identity
\begin{equation}
\label{energy}
\del_t \left ( \frac{1}{2} |v|^2 + W(F) \right ) - \del_\alpha \left (  v_i \frac{\del W}{\del F_{i \alpha}}(F)  \right )
=0
\end{equation}
expressing the conservation of mechanical energy. In the parlance of conservation law $\eta = \frac{1}{2} |v|^2 + W(F)$
is an entropy with associated entropy flux $q_\alpha = -  v_i \frac{\del W}{\del F_{i \alpha}}$. The issue in the background is
that, due to the lack of convexity of  $W$,  the entropy $\eta$ is not convex.
On the other hand the standard theory of systems of conservation laws is intricately connected to convexity and is
inapplicable in a direct way. Hyperbolicity of  \eqref{elassys} is equivalent to the condition on $W$:
\begin{align*}
\frac{\del^2 W}{\del F_{i \alpha} \del F_{j \beta}} (F) \; \xi_i \xi_j \nu_\alpha \nu_\beta > 0
\quad \forall \xi \ne 0 \, , \, \nu \in \cS^{d-1}
\end{align*}
called  {\em Legendre-Hadamard conditions} and associated to positive-definiteness of the acoustic tensor.

A lot of effort has been devoted into
finding suitable notions to replace the assumption of convexity, starting from works on the calculus of variations associated to
the minimization of potential energy. The following notions have emerged:
$W(F) \; \mbox{ is called {\it  (strictly) rank-1 convex } }$ if
$$
W( \theta F + (1 -\theta)G ) <  \theta W(F) + (1-\theta) W(G) \quad 0 < \theta < 1 , \; \forall F, G \quad \mbox{with  $F - G = \xi \otimes \nu$}.
$$
Strict rank-1 convexity is equivalent to the Legendre-Hadamard conditions for $W$ sufficiently smooth.
$W(F) $ is called {\it polyconvex} if it is expressed in the form
$$
W(F) = g ( F, \cof F , \dt F ) = g \circ \Phi (F) \qquad \mbox{with $g( \Xi)$ convex}.
$$
The latter definition is intimately connected to the notion of null-Lagrangean which is the analog of constant
functions in the calculus of variations. An integrand $\Phi (F)$ is a {\it  null-Lagrangean} iff
$$
\int_{\Omega} \Phi (\nabla y + \nabla \phi ) \, dx = \int_{\Omega} \Phi (\nabla y ) dx
\quad \forall \; y \in W^{1,p} \, , \; \phi \in C^\infty_c \, .
$$

The following equivalences hold:
\begin{itemize}
\item[]  \qquad $\Phi (F)$ is a {\it  null-Lagrangean}
\begin{align*}
& \Longleftrightarrow  \quad
\int_{\Omega} \Phi (F + \nabla \phi ) \, dx = \Phi (F ) \, |\Omega|
\quad \forall \; F \in M^{d \times d} \, , \; \phi \in C^\infty_c
\\
& \Longleftrightarrow  \quad
\Phi(F) \; \mbox{is rank-1 affine}
\\
& \Longleftrightarrow  \quad \Phi(F) = \alpha \,  F + \beta \, \cof F + c \, \dt F
\\
& \Longleftrightarrow  \quad
\del_\alpha \left (   \frac{\del \Phi }{\del F_{i \alpha}} (\nabla y) \right )  = 0  \quad \mbox{ in $\cD'$}
\end{align*}
\end{itemize}
The reader is referred  to \cite{ball77,BCO81}  for precise statements and their proofs,
The introduction of null-Lagrangeans into the theory of elasticity is due to Ericksen \cite{ericksen62} and Edelen \cite{edelen62}
and their role in the mathematical theory of elasticity is pointed out by  Ball \cite{ball77}.
Their importance arises from the fact that null-Lagrangeans are weakly continuous in $W^{1,p}$ \cite{BCO81},
and from their role in the existence theory of equilibrium elasticity \cite{ball77,BCO81}
and in obtaining  symmetrizable extensions for  polyconvex  elastodynamics  \cite{qin,DST01,wagner09}.

Here, we explore their role in the problem of cavity formation.

\section{Radial motions in 3-d elasticity}
\label{sec-radel}

A motion of the form $y(x,t) =w(R,t) \frac{x}{R}$,  where $R=|x|$, $x \in\RR^d$, is called radial. For an elastic material to support
radial motions it  must be isotropic and as always frame indifferent,
that is the stored energy satisfies
\begin{align*}
W(Q F) = W(F)  = W(F Q)  \quad  \forall Q \in SO(d) \, .
\end{align*}
The class of such stored energies has been characterized to be of the form $W(F) = \Phi (v_1, v_2, ..., v_d)$,
where $\Phi : \RR^d_{++} \rightarrow \RR$ is a symmetric function of the  eigenvalues $v_i$ of $+ \sqrt{F^T F}$
the so called principal stretches (see  \cite{TN65}).

For reasons that will be explained later we will restrict attention here to the case of three space
dimensions $d=3$ and  the special example of stored energy
\begin{equation}
\Phi (v_1, v_2, v_3 ) = \frac{1}{2} ( v_1^2 + v_2^2 + v_3^2 ) + h ( v_1 v_2 v_3 ) \, ,
\tag {H$_1$}
\end{equation}
with $h : \R_+ \to \R_+$ being a  $C^3$ convex function that satisfies $h(\delta) \to +\infty$ as $\delta\to 0+$.
This stored energy is of the polyconvex class, and the behavior as $\delta\to 0+$ is placed to avoid that a finite volume is compressed
down to zero and guarantee that the solutions avoid infinite compression and can thus be interpreted as elastic motions.

\subsection{The equations of radial elasticity.}
For radial motions  $y(x,t) =w(R,t) \frac{x}{R}$,  with $R=|x|$, $x \in\RR^d$, the deformation gradient
is of the form
$$
F = w_R \hat x \otimes \hat x +  \sum_{j=2}^d \frac{w}{R} \hat x_j^\perp \otimes \hat x_j^\perp \, ,
$$
where $\hat x = \frac{x}{R}$, $\hat x_j^\perp$ the perpendicular unit vectors,
and the principal stretches are $w_R$ in the radial direction and $\frac{w}{R}$ of multiplicity $d-1$ in
the orthogonal directions.
The Piola-Kirchhoff stress tensor is computed as
$$
S = \frac{\del W}{\del F} = \frac{\del \Phi}{\del v_1} \hat x \otimes \hat x + \sum_{j=2}^d  \frac{\del \Phi}{\del v_j}
\hat x_j^\perp \otimes \hat x_j^\perp \, .
$$
To represent a physically realizable motion, we impose that  $\det F = w_R(w/R)^{d-1}  >0$ where $F = \nabla y$.
This is equivalent to $w_R > 0$ a condition that also suffices  to avoid  interpenetration of matter for radial motions.

A cumbersome but straightforward computation shows that $w( R,t)$ satisfies the partial differential
equation of second-order
\begin{equation}
\label{radialelas}
\begin{aligned}
       w_{tt}
        &=\frac{1}{R^{d-1}}  \, \del_R \left ( R^{d-1} \frac{\del \Phi}{\del v_1} \big ( w_R,\frac{w}{R}, ... , \frac{w}{R} \big ) \right )
         -\frac{d-1}{R}  \frac{\del \Phi}{\del v_2}   \big (w_R,\frac{w}{R}, ..., \frac{w}{R} \big )  \, .
\end{aligned}
\end{equation}

It is instructive to give a short alternative derivation of \eqref{radialelas}. Consider the action functional $\mathcal{L} = K -P$,
where $K$ is the kinetic and $P$ the potential energy:
$$
\mathcal{L} [w] =  K - P = \int_0^T \int_0^1 \frac{1}{2} w_t^2 R^{d-1} -  \Phi ( w_R , \frac{w}{R} , ... , \frac{w}{R} ) R^{d-1} dR dt\, .
$$
The critical points of the functional $\mathcal{L}$ are computed by setting to zero the variational derivative
 $\frac{d}{d\delta} \Big |_{\delta = 0} \mathcal{L} (w + \delta \psi) = 0 $.
 This leads to the  Euler-Lagrange equations
$$
\int_0^T \int_0^1  R^{d-1} w_t  \psi_t -  R^{d-1} \frac{\del \Phi}{\del v_1} \psi_R
- (d-1) \frac{\del \Phi}{\del v_2} \frac{\psi}{R} R^{d-1} \;  dR dt  = 0
$$
which are precisely the weak-form of \eqref{radialelas}.

\subsection{Radial elasticity features}

Upon introducing the velocity, radial strain, and transversal strain, respectively
$$
v = w_t \quad a = w_R \quad b = \frac{w}{R}
$$
\eqref{radialelas} is expressed as the first order system
\begin{equation}
\label{radialelassys}
\begin{aligned}
       \del_t v
        &=\frac{1}{R^{d-1}}  \, \del_R \left ( R^{d-1} \frac{\del \Phi}{\del v_1} \big ( a , b, ...,  b \big ) \right )
         -\frac{d-1}{R}  \frac{\del \Phi}{\del v_2}  \big (a, b , ...,  b \big )
\\
       \del_t a &= \del_R v
\\
        \del_t b &= \frac{1}{R} v
\end{aligned}
\end{equation}
subject to the constraint $\del_R (Rb) = a$, which is again an  {\it involution} propagating from the initial data and causes no
problems as constraint. For $\Phi_{11} : = \frac{\del^2 \Phi}{\del v_1^2} > 0$, the system \eqref{radialelassys} is hyperbolic
with wave speeds  $\lambda_{\pm} = \pm \sqrt{\Phi_{11}(a,b, ..., b)}$,  $\lambda_0 = 0$.

A comparison  with the system of elasticity in one-space dimension \eqref{elasoned}  indicates
the following differences
\begin{itemize}
\item[(i)]  The system \eqref{radialelassys} has lower order terms and a geometric singularity at the origin.

\item[(ii)]   As contrasted to the case of one-dimensional elasticity the  wave speeds $\lambda_{\pm}$ depend on the
lower order terms.

\item[(iii)]   The system has an involution.

\item[(iv)]   Following the usual theory of conservation laws a weak entropy solution will be defined to satisfy
the entropy inequality
\begin{equation}
\label{entrsol}
\del_t \left (  \frac{1}{2} v^2 + \Phi (a, b, ... , b)   \right  ) R^{d-1}  + \del_R \left ( R^{d-1}  v  \frac{\del \Phi}{\del v_1} (a,b, ..., b) \right ) \le 0
\end{equation}
expressing the dissipation of mechanical energy.  The assumption \eqref{hyp1} implies that,
in a deviation to the usual theory,  the entropy here is merely polyconvex.
\end{itemize}

Little is known at present concerning the existence of weak solutions for the radial elasticity
system \eqref{radialelassys}. One additional difficulty is that solutions of  \eqref{radialelassys} have to be constructed
to satisfy the constraint $\det F >0$ so as to be interpreted as mechanical motions. This constraint is difficult to
preserve even in approximating schemes, with the  exception of the variational
approximation scheme in  \cite{MT12} which preserves the positivity of Jacobians and produces iterates that
decrease the mechanical energy.


\section{ The problem of cavitation}
\label{sec-cavi}

\begin{figure}
    \centering\includegraphics[height=3cm,width=6.5cm]{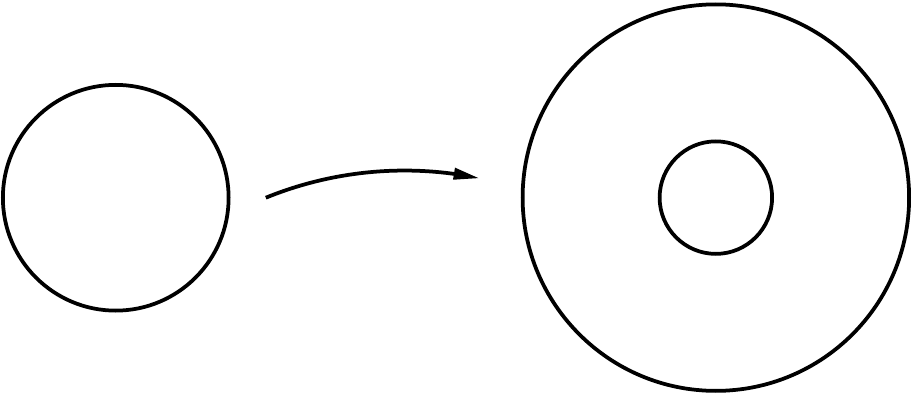}
    \caption{A cavitating solution}
    \label{fig2}
  \end{figure}

The objective is to consider the equations of elasticity
$$
\frac{\del^2 y }{\del t^2} = \div \frac{\del W}{\del F} (\nabla y)
$$
and to examine under what conditions the equations admit solutions with cavities.
As already noted such solutions  are radial,  $y (x,t) = w(R,t) \frac{x}{R}$ with $R = |x|$, and there will be a cavity provided $w(0,t) > 0$,
see Figure \ref{fig2} for a depiction of such a solution.

The natural energy bound for the elasticity system is
$$\int \frac{1}{2} |y_t|^2  +  |\nabla y|^p   dx < \infty $$
with the exponent $p$ determined by the growth of the stored energy $W$.
Accordingly, if one is interested to explore cavitating solutions with finite energy, the Sobolev
embedding theorem dictates to consider growth exponents $p < d$, where $d$ is the space dimension.
Cavitating solutions lie outside the realm of the usual continuum modeling, and one might question if
continuum mechanics can account for their presence. Nevertheless, Ball \cite{ball82} in a groundbreaking work examined whether
the equations of equilibrium elasticity admit radial solutions with cavities.
He noted that
$$
y(\cdot, t)  \in W^{1,p}_{loc} \Longleftrightarrow \left \{
\begin{aligned}
w(\cdot , t) \; \mbox{ absolutely continuous on $(0,1)$}
\\
\int_0^1 \Big ( |w_R|^p + \frac{|w|^p}{R^p} \Big ) R^{d-1} dR < \infty
\end{aligned}
\right .
$$
and that for $y \in W^{1,1}_{loc}$ and for $d \ge 2$ there is no delta-mass at the origin associated to the cavity, i.e.
$$
\begin{aligned}
&\nabla y  = \frac{w}{R} I  + ( w_R  - \frac{w}{R} )  \hat x \otimes \hat x
\quad \mbox{ in $\mathcal{D}'$ and a.e.}
\end{aligned}
$$
This situation should be contrasted with the case $d=1$ where a jump discontinuity induces a delta
mass at the origin.

\subsection{Cavitation in equilibrium elasticity}
Ball \cite{ball82} studied the cavitation problem from two perspectives. Using direct methods of the calculus of variations,
he studies the minimization problem
$$
\min \int_0^1 \Phi ( w'(R), \frac{w}{R}, ... , \frac{w}{R} ) \, R^{d-1} dR
$$
over the set of admissible functions $\mathcal{A}_\lambda$.
He also considered the
associated Euler-Lagrange equations
\begin{equation*}
\begin{aligned}
       \, \del_R \left ( R^{d-1} \frac{\del \Phi}{\del v_1} \right )
         - (d-1)  R^{d-2} \frac{\del \Phi}{\del v_2}   &= 0
         \\
      w(1) &= \lambda
\end{aligned}
\end{equation*}
and carried out a bifurcation analysis for this system. The analysis in \cite{ball82} shows:

\begin{itemize}
\item There is an absolute minimum in the admissible set
$$
\mathcal{A}_\lambda = \{ w \in W^{1,1} , w(0 ) \ge 0, w'(R) > 0 , w(1) = \lambda \} \, .
$$

\item The minimizer $w_\lambda$ satisfies the Euler-Lagrange equations for $R > 0$ and
\begin{itemize}
\item[a.]  There is a critical value $\lambda_{cr}$ such that for $\lambda > \lambda_{cr}$ there is a unique
minimizer $w_\lambda$ with $w_\lambda (0) > 0$.
\item[b.] The trivial solution $w = \lambda R$ is stable (in a suitable sense) if $\lambda \le \lambda_{cr}$ and unstable
if $\lambda >  \lambda_{cr}$.  If $\lambda > \lambda_{cr}$ then $w_\lambda$ is stable.
\end{itemize}
\end{itemize}
The reader is referred to Ball \cite{ball82} for the details (and precise statements) of these results and to
Sivaloganathan and Spector \cite{SS03} for an account of later results on the problem of cavitation in
elastostatics. The problem of cavitation has also been studied in a context of incompressible elasticity \cite{ball82}
which lies outside the realm of our discussion.

\subsection{Dynamic radial elasticity - Cavitation}

Consider next the problem of dynamic cavitation in compressible, isotropic elastic materials.
Such materials have stored energies $W(F) = \Phi (v_1, v_2, v_3)$ and support an {\it ansatz} of radial motions
with $w(R,t)$ satisfying the equations of radial elasticity \eqref{radialelas}. The  latter admit the special solution
$w(R,t) = \lambda R$ which corresponds to a homogeneous elastic deformation of stretching $\lambda$ and is a
special solution emanating from initial data $w_0 (R) = \lambda R$.

Pericak-Spector and Spector in a remarkable work \cite{ps88,ps97} postulated a {\it self-similar ansatz} of solutions
\begin{equation}
\label{sscav}
y(x,t) = w(R,t) \frac{x}{R} = t \varphi \big ( \frac{R}{t} \big ) \frac{x}{R}
\end{equation}
and constructed a second self-similar solution for the dynamic elasticity system associated to cavity formation $\varphi(0) > 0$
provided the stretching $\lambda$ is bigger than some critical value, $\lambda > \lambda_{cr}$.

\begin{theorem}
[Pericak-Spector and Spector \cite{ps88}]
For a stored energy function
\begin{equation}
 \label{hyp1}
\Phi (v_1, v_2, v_3 ) = \frac{1}{2} ( v_1^2 + v_2^2 + v_3^2 ) + h ( v_1 v_2 v_3 ) \, ,
\tag {H$_1$}
\end{equation}
with $h: \R_+ \longrightarrow \R_+$ satisfying the hypotheses
\begin{equation}
 \label{hyp2}
 h''>0 \, ,  \quad  h'''<0 \, \quad  \lim_{v \rightarrow 0} h(v) = \lim_{v \rightarrow \infty} h(v) = \infty  \, ,
 \tag {H$_2$}
\end{equation}
and for $\lambda > \lambda_{cr}$ there exist cavitating solutions \eqref{sscav} for dimension $d = 3$ satisfying
\begin{equation*}
\begin{aligned}
       w_{tt}
        &= \frac{1}{R^2} \del_R \left ( R^{2} \frac{\del \Phi}{\del v_1} \big ( w_R,\frac{w}{R},\frac{w}{R} \big ) \right )
         -\frac{1}{R} \Big ( \frac{\del \Phi}{\del v_2} + \frac{\del \Phi}{\del v_3} \Big )
          \big (w_R,\frac{w}{R},\frac{w}{R} \big ) \; ,    \; R > 0 \, , \; t>0
  \\
       w(R,t) &= \lambda R  \quad  \mbox{for  $R > \bar{r }t $ for some $\bar{ r } > 0$}
\end{aligned}
\end{equation*}
and $w(0, R) = \lambda R$.
\end{theorem}

An initial version of this theorem was proved in  \cite{ps88} for \eqref{hyp1} and $d \ge 3$ and it was extended
in \cite{ps97} to a far more general class of polyconvex constitutive functions $\Phi$.
A variant also holds for dimension $d=2$ but the growth of $\Phi$
is there restricted to be slightly superlinear.
Additional regularity properties for the cavitating solutions and a bifurcation study is provided in \cite{MT13}.

We give an outline of the main ideas in the construction and the properties of the constructed solutions. Our exposition
follows \cite{ps88} and \cite{MT13} and we refer there for the details.
The starting point is the self-similar {\it ansatz}
\begin{equation}
\label{ssansatz}
\begin{aligned}
w(R,t) &= t \varphi \big ( \frac{R}{t} \big ) \qquad s = \frac{R}{t} \, .
\end{aligned}
\end{equation}
Note that $w(0,t) = t \varphi(0)$ and thus $\varphi (0)$ stands for the velocity of the cavity. Introducing
\eqref{ssansatz} to \eqref{radialelas} leads to the ordinary differential equation
\begin{equation}
\label{mainsseq}
\begin{aligned}
(s^2 - \Phi_{11} ) \ddot{\varphi} = \frac{2}{s} ( \dot \varphi - \frac{\varphi}{s} )
 \underbrace{ \left [ \Phi_{12} + \frac{ \Phi_1 - \Phi_2}{\dot \varphi  - \frac{\varphi}{s}}
         \right ] }_{P ( \dot \varphi , \frac{\varphi}{s} , \frac{\varphi}{s} ) \ge 1 } .
\end{aligned}
\end{equation}
An important property is that due to \eqref{hyp1} the function $P$ defined above satisfies $P \ge 1$.

It is instructive to rewrite \eqref{mainsseq} as a first order system for the variables
$$
a = \dot \varphi \, , \quad b = \frac{\varphi}{s}
$$
It takes the form
\begin{equation}
\label{mainsys}
\begin{aligned}
\big ( s^2 - \Phi_{11} (a,b,b) \big ) \dot a &= \frac{2}{s} ( a -b ) P(a,b,b)
\\
\dot b &= \frac{1}{s} ( a- b)
\end{aligned}
\end{equation}
\begin{figure}
\centering
 (a) \includegraphics[height=4cm,width=6cm]{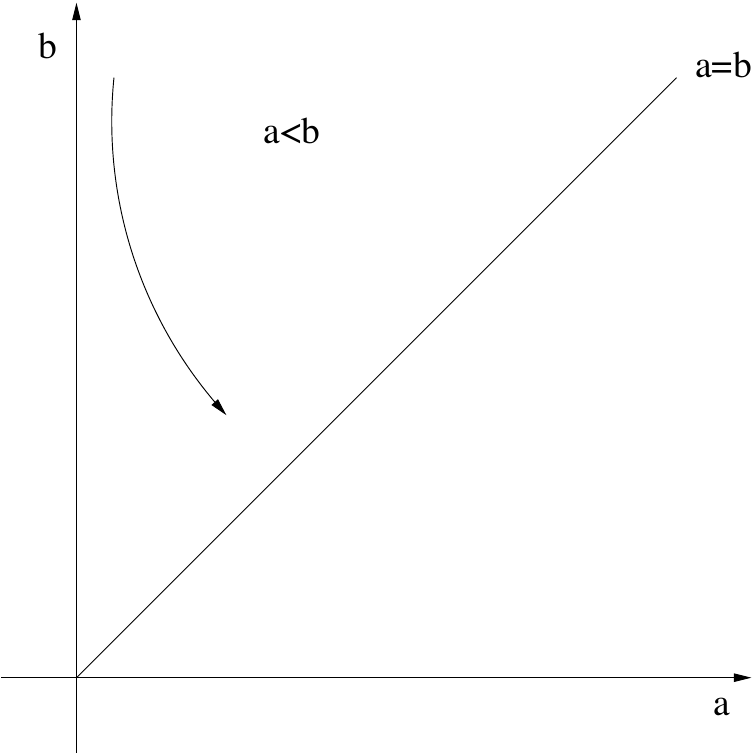} \hfil
  (b) \includegraphics[height=4cm,width=6cm]{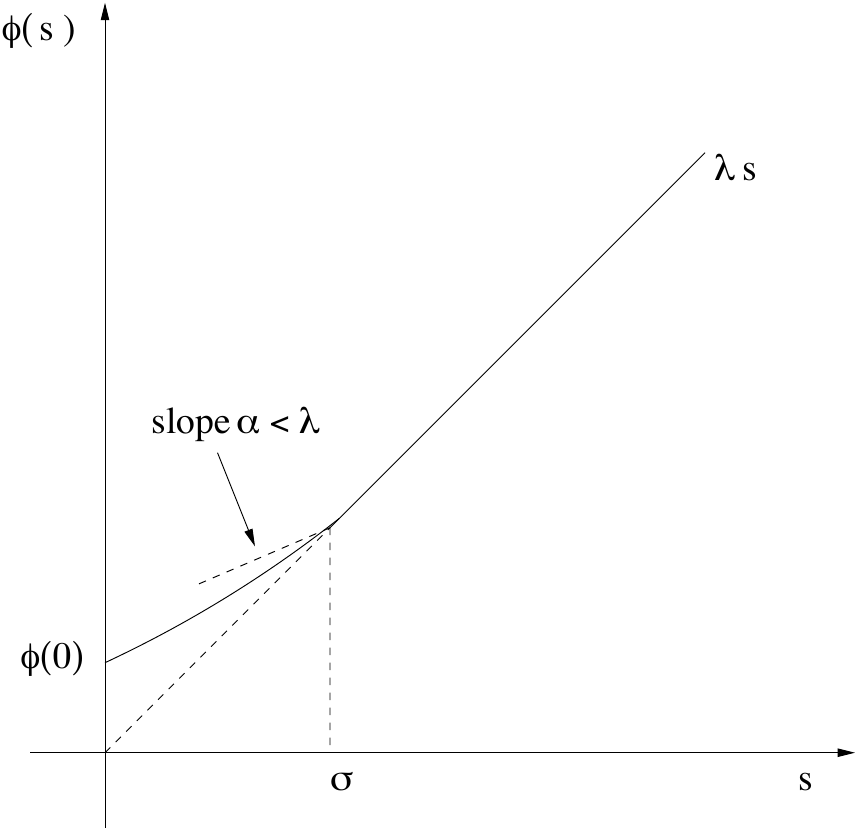}
    \caption{(a) A typical trajectory of  \eqref{mainsys}; (b) The cavity solution}
    \label{fig3}
 \end{figure}

A simple examination of the system \eqref{mainsys} reveals the following properties:
\begin{itemize}
\item  The system is singular with a free boundary type of singularity at the sonic line $s^2 - \Phi_{11} = 0$.

\item The system has one class of special solutions, the states $a = b$, corresponding to a homogeneous deformation.

\item A second  class of solutions consists of jump discontinuities. These satisfy the Rankine-Hugoniot conditions
$$
\begin{aligned}
\sigma^2 [a] - [\Phi_1] &= 0
\\
[b] &= 0
\end{aligned}
$$
and give rise to shocks. Note that along a shock, $b$ is continuous but $a$ jumps. For functions of class \eqref{hyp1}, hypothesis
\eqref{hyp2} implies that $\Phi_{111} < 0$ and it can be worked out that the Lax shock admissibility criterion implies that
the admissible shocks satisfy  $b_- = b_+$ and  $a_- < a_+ = \lambda$.

\item The main item is to study the third class of solutions to \eqref{mainsys}, that is  continuous solutions.
The analysis is based on a monotonicity property that can be read directly from \eqref{mainsys}. So long as
the solution does not cross the equilibrium diagonal $a=b$ or the sonic singularity, solutions are monotone functions, see Figure \ref{fig3}(a).
\end{itemize}

\noindent
The analysis of the continuous solution of \eqref{mainsseq} is based on a de-singularization of the problem. One
may introduce the variables $(\varphi , v)$ where $v = \dot \varphi ( \frac{\varphi}{s} )^2$ and rewrite \eqref{mainsys} in
the form
$$
\begin{aligned}
\dot \varphi &= \hat \Phi (s, \varphi, v)
\\
\dot v &= \hat U (s, \varphi, v)
\end{aligned}
\qquad
\left \{
\begin{aligned}
\varphi(0) &= \varphi_0 > 0
\\
v(0) &= v_0 > 0 \, .
\end{aligned}
\right .
$$
Under \eqref{hyp1} the functions $\hat \Phi$ and $\hat U $ are {\it non-singular} at the origin $s=0$ and the standard existence
theory can be applied to the last system. At the origin there hold the following properties for the Cauchy stress $T(0)$ (see \cite{ps88})
$$
\begin{aligned}
\mbox{either} \quad  &T(0) = 0 \Leftrightarrow v(0) = H   \quad  \mbox{stress free cavity}
\\
\mbox{or} \quad  &T(0) = G (\varphi(0))  \quad  \mbox{cavity with content} .
\end{aligned}
$$
The second condition is analogous to the kinetic relations familiar from the theory of phase transitions.
We refer to \cite{MT13} for the corresponding analysis of the singularity at dimension $d=2$.

The constructed solution has the following properties:
\begin{itemize}
\item[(i)] $a = \dot \varphi \nearrow$ \quad $b = \frac{\varphi}{s} \searrow$ \quad $a - b  \nearrow$.
\item[(ii)]  It can be extended on a maximal interval of existence $(0, T)$
$$
Q = s^2 - \Phi_{11} \to 0 \qquad a-b \to c < 0 \qquad \mbox{as $s \to T$}.
$$
\vskip-0.3cm
\item[(iii)]  If $\Phi_{111} < 0$  the solution is connected to a uniformly deformed state through a single shock
that is admissible via the Lax shock admissibility condition.
\end{itemize}
We refer to \cite{ps88} and \cite{MT13} for the proofs. It should be noted that the constructed cavitating solution
consists of a single precursor shock that connects a convex function $\varphi(s)$ to a uniform deformation, see Figure \ref{fig3}(b).
It is necessary to have a precursor shock, that is, it is not possible for $\varphi(s)$ to connect in a $C^1$-fashion
to a uniformly deformed state via a sonic singularity (see \cite{MT13}).
In summary, $\varphi(s)$ looks like:
\begin{itemize}
\item[(a)]  At $\sigma$ there is a shock.
\item[(b)]  $\varphi(s)$ convex for $s < \sigma$.
\item[(c)]  $\varphi(s) = \lambda \, s $ for $s > \sigma$.
\end{itemize}

The following remarkable property was proved in \cite[Thm 7.2]{ps88}:
The cavitating solution {\it decreases} the mechanical energy, namely, if
$$
E( y, B_\rho)  = \int_{B_\rho}  \frac{1}{2} |y_t|^2 + W (\nabla y) dx
$$
$y_h = \lambda x$ is the homogeneous solution  and  $y_c$ is the solution with the cavity, then
\begin{equation}
\label{cavdec}
\begin{aligned}
E (y_c , B_\rho)  - E (y_h , B_\rho) &= (t \sigma)^3 \frac{4 \pi}{3}
\Big [ \Phi (a_- , \lambda, \lambda) - \Phi (\lambda, \lambda, \lambda)
\\
&\quad+\frac{1}{2} \big ( \Phi_1 (a_- , \lambda, \lambda) + \Phi_1 (\lambda , \lambda, \lambda) \big )
(\lambda - a_-) \Big ]
\\
&< 0 \qquad \mbox{whenever} \quad a_-  := \varphi( \sigma -)  < \lambda
\end{aligned}
\end{equation}
where the right hand side in \eqref{cavdec} corresponds to the dissipation at the precursor outgoing shock and
is strictly negative.

This result states nonuniqueness for entropy weak solutions  (with polyconvex energy)
because of the singularity at the cavity. As already noted in \cite{ps88},
the paradox arises that by opening a cavity the energy of the material decreases,
what induces an autocatalytic mechanism for failure.
From a  perspective of  mechanics the problem appears to be that there is no surface energy cost for opening the cavity
that is accounted for by the weak solution.

\section{ Fracture in 1-d}
\label{sec-frac}

Insight into the non-uniqueness issue can be obtained by studying the equations of one-dimensional
elasticity or the equivalent form of the elasticity system,
\begin{equation}
\label{elasoned}
y_{tt} = \del_x \tau (y_x )
\qquad
\longleftrightarrow
\qquad
\begin{cases}
u_t - v_x = 0 & \\
v_t - \tau (u)_x = 0 & \\
\end{cases}
\, ,
\end{equation}
where $v= y_t$ is the velocity and $u= y_x$ is the strain.
This equation admits the special solution $y_h (x,t) = \lambda x $ corresponding to a homogeneous
deformation with strain $u_h = \lambda$ and velocity $v_h = 0$.
We assume that the stress function $\tau (u)$ satisfies the hypotheses
\begin{equation}
\label{hypo1}
\tau'(u) > 0 \, , \quad \tau''(u) < 0
\tag{$a_1$}
\end{equation}
\begin{equation}
\label{hypo2}
\tau(u) \to - \infty \quad \mbox{as $u \to 0$} \quad \mbox{and} \quad \int_1^u \tau(s) ds \to + \infty \quad \mbox{as $u \to 0.$}
\tag{$a_2$}
\end{equation}
Under \eqref{hypo1} the wave speeds $\lambda_{1,2} (u) = \pm \sqrt{ \tau'(u) }$ are real and \eqref{elasoned} is hyperbolic.
The hypothesis $\tau''(u) < 0$ is appropriate for an elastic material exhibiting softening elastic response
and plays an important role in the forthcoming analysis.
The hypothesis \eqref{hypo2} is applicable in the case of longitudinal motions and
is placed to exclude that a finite
volume is compressed down to zero. In the sequel we will consider only tensile deformations
and this hypothesis will not play any significant role.
 (In the case of shearing motions $\tau (u)$ is defined for $u \in \R$ and \eqref{hypo2} is removed).

Motivated from the problem of cavitation we introduce the {\it ansatz}
\begin{equation}
\label{ssoned}
y(x,t) =   t \varphi \Big ( \frac{ |x| }{t} \Big ) \frac{x}{|x|} = t Y ( \frac{x}{t} )
\end{equation}
with $Y(0) > 0$. This ansatz is similar to the one used for the solution of the Riemann problem for \eqref{elasoned}
except at the origin $x=0$. A calculation shows that
$$
\del_x y = Y' (\xi) + 2 Y(0) \delta_{\xi = 0}
$$
and such an ansatz could conceivably provide a cavitating (in fact fracturing) solution, except that at the origin
there is a delta mass and a suitable interpretation to the term $\tau (\del_x y)$ should be supplied. Of course this is
the main difference between the case $d=1$ and the dimensions $d \ge 2$ where a delta mass in the origin does not appear.

In \cite{GT13}, the following function
\begin{equation}
\label{singsol}
y(x,t) =  t Y ( \frac{x}{t} ) \qquad Y(\xi) =
\begin{cases}
\quad \lambda \xi & \xi < - \sigma \\
-Y(0) + \alpha \xi & -\sigma < \xi < 0 \\
\; \; Y(0) + \alpha \xi  &  0 < \xi < \sigma   \\
\quad \lambda \xi & \sigma < \xi
\end{cases}
\end{equation}
is tested as a candidate for solution of \eqref{elasoned}. The constants $\alpha, \, \lambda$ and $Y(0)$ are selected to sastisfy
the Rankine-Hugoniot jump conditions at the two outgoing shocks at $\xi = \pm \sigma$
$$
\begin{aligned}
Y(0) &= \sigma (\lambda - \alpha)
\\
\sigma^2 &= \frac{\tau(\lambda) - \tau(\alpha)}{\lambda - \alpha}
\end{aligned}
$$
One easily checks that
$$
\mbox{ $\tau''(u) < 0$ and $\alpha < \lambda$ $\Longrightarrow$ both shocks are Lax-admissible. }
$$

\begin{figure}
    \centering\includegraphics[height=5cm,width=8cm]{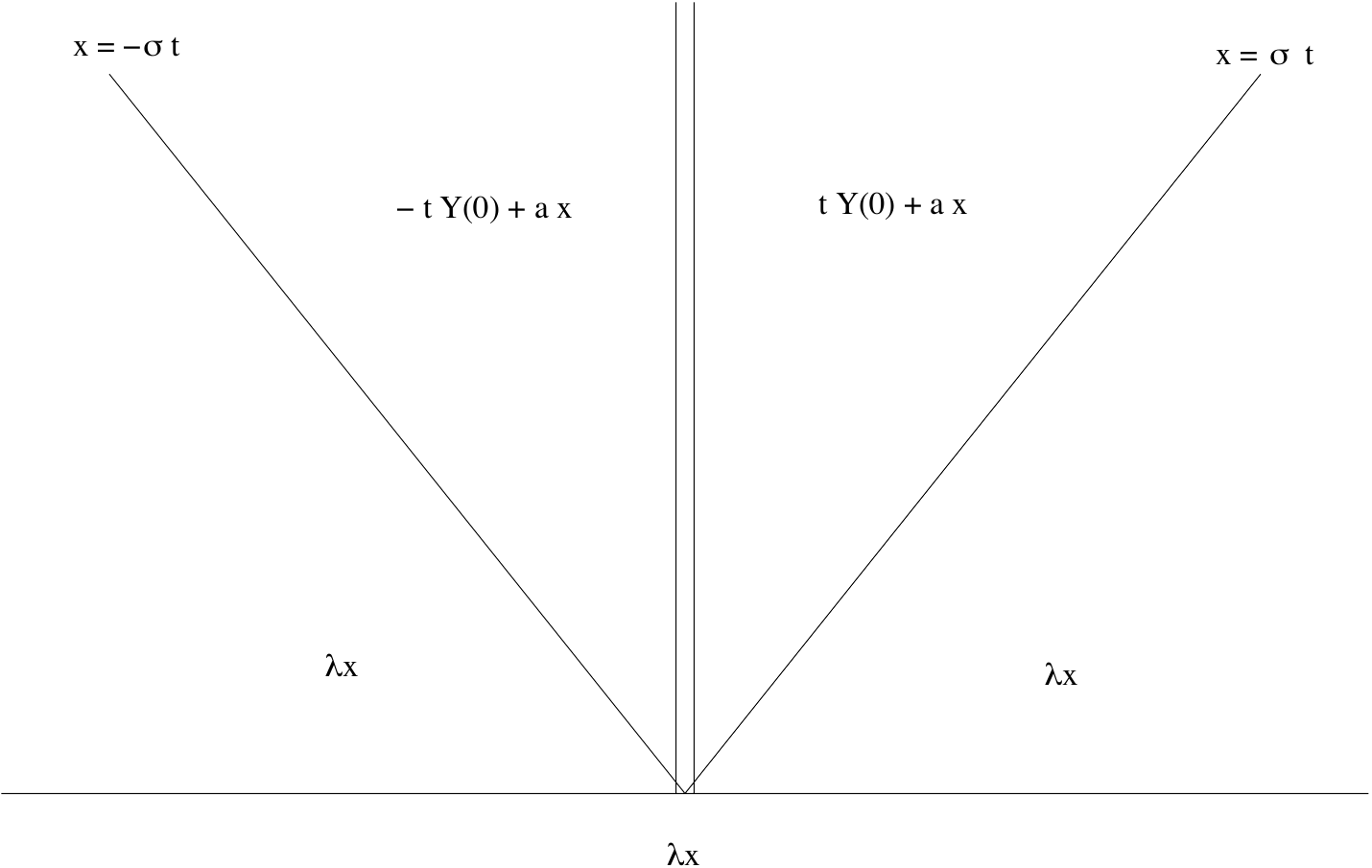}
    \caption{A fracturing motion in Lagrangian coordinates}
    \label{fig5}
  \end{figure}

\begin{figure}
    \centering\includegraphics[height=5cm,width=8cm]{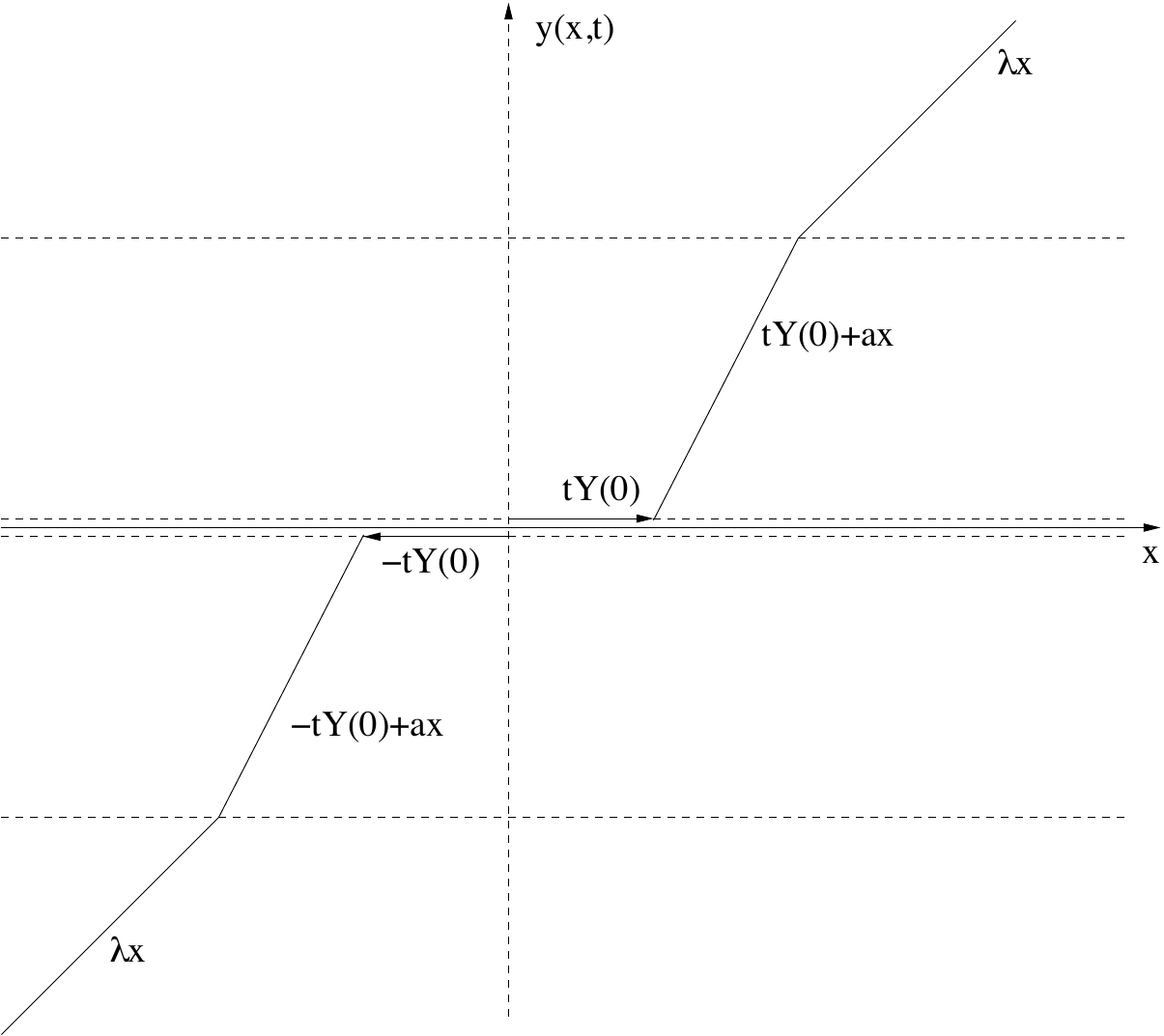}
    \caption{A shear band in Eulerian coordinates}
    \label{fig6}
  \end{figure}

The question arises if \eqref{singsol} can be given a suitable interpretation as a solution of \eqref{elasoned}.
Note that if the equation \eqref{elasoned} is interpreted as longitudinal motions then such a solution would
correspond to {\it fracture}, see Figure \ref{fig5}.
If \eqref{elasoned} is interpreted as shear motions then \eqref{singsol} corresponds to a {\it shear band}, see Figure \ref{fig6}.
The difference among the two cases is that the material splits apart in the former case, whereas it remains in contact in the latter.
A more precise statement of the underlying question is whether the model \eqref{elasoned} with the hyperbolic and strain softening
stress-strain relation \eqref{hypo1} may support such motions.

\section{ Slic-solutions}
\label{sec-slic}

There is a class of problems in material science where structures with discontinuous displacement fields
emerge via a material instability mechanism. Typical examples are development of cracks in fracture,
cavitation in elastic response, or formation of shear bands in plastic deformations. Such problems lie at
the limits of applicability of continuum modeling and the usual concept of weak solutions is in
any case inadequate to describe these motions. Nevertheless, as the material transitions
from a regime where continuum modeling is applicable to a regime that the model has to be modified (or perhaps
atomistic modeling has to be employed), it is expected that at the interface both types of modeling have to apply
 in an intermediate regime. It is further expected that such structures should appear in a small parameter limit of more complex models
that incorporate "higher-order physics", and that their appearance occurs in a stable way.

The concept of  {\it singular limiting induced from continuum solution} (in short {\it slic}-solution) \cite{GT13}
is an attempt to give meaning to such discontinuous solutions. Roughly speaking a discontinuous solution of \eqref{elasoned}
will be a slic-solution if it can be obtained as the limit of approximate smooth solutions that are an averaging of $y$.
Due to the self-similar form of \eqref{ssoned} averaging in space entails averaging in time.
In \cite{GT13}, the following definition is introduced and used to test the example \eqref{singsol}.

\begin{definition}
\label{defslicconvx}
Let $y \in L^\infty_{loc} \big ( (-\infty ,\infty) \, ; \, L^1_{loc} (\R) \big )$
 satisfy for some $\eps  > 0$
the monotonicity condition:
\begin{equation}
\label{hypomc}
y(x_1,t)-y(x_2,t)>  \eps (x_1-x_2)  \quad \mbox{for $x_1,x_2,t \in \R$ with $x_1>x_2$.}
\tag{$mc$}
\end{equation}
For  $\phi$ a mollifier,
 $\phi \in C^\infty_c (\R)$, $\phi \ge 0$,  $\supp \phi \subset  B_1$ (the ball of radius 1),
$\int \phi = 1$, we let $\phi_n = n  \phi  ( n x   )$ and  define the averaged function
\begin{equation}
\label{convx}
y^n (x,t) = \phi_n \underset{x}{\star}  y = \int \phi_n (x-z) y(z,t) dz.
\end{equation}
The function $y$ is called a singular limiting induced from continuum (slic-) solution of \eqref{elasoned} provided
for any symmetric mollifier
\begin{equation}
\label{eqdefslic}
\del_{tt} y^n - \del_x \tau (\del_x y^n ) =: f^n \to 0 \quad \mbox{ in $\cD'$}
 \end{equation}
as $n \to \infty$.
\end{definition}
The example \eqref{singsol} is extended to negative values of $t$ by setting $y = \lambda x$ for $t < 0$. The resulting
function is given in explicit form by
\begin{equation*}
y(x,t) =
 \begin{cases}
\lambda x  \charf_{x < -\sigma t}
+ ( - t Y(0) + \alpha x)  \charf_{- \sigma t < x < 0}
+ ( t Y(0) + \alpha x)  \charf_{0 < x < \sigma t}
            + \lambda x  \charf_{\sigma t < x}
  & t > 0
  \\
  \lambda x & t < 0 \, .
  \end{cases}
\end{equation*}
Using the explicit form, we calculate the averagings in \eqref{convx}.
A lengthy computation in \cite{GT13}
shows that for  $ \psi (x,t)  \in C_c^\infty (\R \times \R)$
$$
\begin{aligned}
 \int_\R \int_\R y^n \psi_{tt} + \tau(  y^n_x ) \psi_x \, dx dt  &=
\int_0^\infty \int_{-\frac{1}{n}}^{\frac{1}{n}}  \tau \big ( \alpha + 2 \phi_n (x) t Y(0) \big ) \psi_x  dx dt  \; +  \; o(1)
\\
&\to 2 t Y(0) L \int_0^\infty \psi_x (0,t) dt
\end{aligned}
$$
where $L = \lim_{u \to \infty} \frac{\tau (u)}{u} $. One thus  concludes:
\begin{itemize}
\item[] If  $L > 0$ then $y(x,t)$ is not a slic-solution.
\item[] If $L = 0$ then $y(x,t)$ is a slic-solution.
\end{itemize}
Consequently, a meaning to \eqref{singsol} as a slic-solution can only be given for sublinear growth of $\tau (u)$.

Next, we consider the energy balance.  Let  $B = (-r, r)$ contain the entire wave fan of the approximate solution
\eqref{convx} at time $t$.  The velocity at the boundary of such a domain vanishes, $v^n \big |_{\del B} = 0$ and the
total energy of the wave fan
\begin{equation}
\label{toten}
E [ y^n ; B  ] = \int_B \frac{1}{2} (\del_t y^n)^2 + W( \del_x y^n)  \, dx
\end{equation}
evolves according to the energy balance equation
$$
\frac{d}{dt} \int_B \frac{1}{2}(\del_t y^n)^2 + W( \del_x y^n)  dx = \int_B f^n v^n \,dx \, .
$$
A calculation again shows
$$
\begin{aligned}
\frac{d}{dt} \int_B \frac{1}{2} (\del_t y^n)^2 &+ W (\del_x y^n ) dx = \int_B f^n v^n dx
\\
&= Y(0)^2 \sigma - 2 \sigma ( W(\alpha) - W(\lambda) ) + 2 \int_0^{\frac{1}{n}}
\tau \big ( \alpha + 2 \phi_n(x) t Y(0) \big ) 2 Y(0) \phi_n (x) dx
\\
&\to \mu_{-\sigma} + \mu_{\sigma} + p_c =: T
\end{aligned}
$$
where the total energy production $T$ is split into the energy dissipation of the two shocks
$\mu_{\pm \sigma}$ and the contribution to the energy by the surface energy of the cavity $p_c$.
One easily also checks that
$$
\begin{aligned}
T &= \mu_{-\sigma} + \mu_\sigma + 2 (\tau_\infty - \tau(\alpha)) Y(0)
\\
&= \sigma Y(0)^2 - 2 \sigma (W(\lambda) - W(\alpha)) + 2 \tau_\infty Y(0)
\\
&= \sigma Y(0)^2 + 2 Y(0) \Big ( \tau_\infty - \frac{W(\lambda) - W(\alpha)}{\lambda - \alpha} \Big )
> 0
\end{aligned}
$$
so that if  $\tau_\infty = \infty$ then $T = + \infty$ while if $\tau_\infty < \infty$  then $0 < T < \infty$.

\section{Cavitation in 3-d}
\label{sec-sliccav}

Next, we consider solutions of the isotropic elasticity equations \eqref{radialelas}, with stored energy
satisfying \eqref{hyp1} and \eqref{hyp2}, of the form
$$
y(x,t) =  w(|x|,t) \frac{x}{|x|} = t  \varphi(s)  \frac{x}{|x|} \, , \quad  s = \frac{|x|}{t}
$$
where $\varphi(s)$ is the cavitating self-similar solution outlined in section \ref{sec-cavi}.
The objective is to examine how the intuition from the example in section \ref{sec-frac}
transfers to the problem of cavitation.

The following extension of the notion of slic-solution -- adapted to the radial case  -- can be introduced.

\begin{definition}
\label{def:slic3d}
Let $y(x,t) = w(|x|,t) \frac{x}{|x|}$ with $w \in L^\infty_{loc}(\R; L_{loc}^1(\R) )$ and $w(\cdot, t)$ monotone increasing
satisfy $y(x,t) = \lambda x$ for $t \le 0$ and for $|x| > \bar{r} t \, , t >0$ or some $\bar r>0$.
The function $y$ is  called a  singular limiting induced from continuum ({\it slic})-solution  of \eqref{elas-m}
if
\begin{equation}
\label{slic-av-3d}
y^n = w^n (|x|,t) \frac{x}{|x|} \, , \quad \mbox{ with \; \; $w^n = \phi_n \underset{R}{\star}  w$} \, ,
\end{equation}
satisfies $\det \nabla y^n \ge \eps_n > 0$ for all $\in \mathbb{N}$, and
$$
\frac{\del^2 y^n }{\del t^2} - \div \frac{\del W}{\del F} (\nabla y^n ) =: f^n \to 0 \quad \mbox{ in $\cD'$}
$$
as $n \to \infty$, for all $\phi \in C_c^\infty(\R)$ positive, symmetric mollifiers with $\phi(0) > 0$.
\end{definition}

The convolution $w^n = \phi_n \underset{R}{\star}  w$ is defined by first extending $w(R,t)$ antisymmetrically in $R$
and the definition takes into account the layer structure at the cavity.
In agreement with the definition of slic-solution (and of the energy in one space dimension) we define the energy of a multidimensional slic-solution.

\begin{definition}\label{def:energy3d}
Let $w \in  W^{1,\infty}_{loc}(\R; L_{loc}^1( \R ) )$, The energy of a slic-solution
in a domain $B \subset \R^d$ containing the entire wave fan and for a.e. $t \in \R$ is defined as
$$
E[ y ,B] (t) :=  \lim_{n \to \infty}   \int_B \frac{1}{2} |y^n_t(x , t )|^2 + W(\nabla y^n(x,t)) dx
$$
with $y^n$ defined by \eqref{slic-av-3d}.
\end{definition}

The following is proved in \cite[Thm 3.8, Prop 3.10]{GT13} based on detailed estimations of $w^n = \phi_n \underset{R}{\star}  w$
on the entire wave fan but mainly around the cavity:

\begin{theorem} (i) The weak solution constructed by Pericak-Spector and Spector is a slic-solution  if and only if
$$
\lim_{u \to \infty} \frac{ h'(u^3 ) }{ u } =0 \, .
$$
(ii)  If $L = \lim_{u \to \infty} \frac{h(u)}{u}$ then
$$
\begin{aligned}
\lim_{n\to \infty} \int_B \frac{1}{2} &|v^n|^2 + W (\nabla y^n) - W(\lambda) dx
\\
&= \underbrace{\int_B  \big [ \frac{1}{2} |v|^2 + W (\nabla y) - W(\lambda) \big ]  dx}_{ < 0}
+ \underbrace{  \big ( t \varphi(0) \big )^3 \,  \frac{4 \pi }{3}   L}_{> 0}  = :  P_{wf} > 0  \, .
\end{aligned}
$$
\end{theorem}
In particular, (ii) shows that for a slic-solution the energy of the solution with cavity  is larger than the energy of the uniform deformation.

\medskip
\noindent
{\bf A rough explanation of the discrepancy}.
We refer to \cite{GT13} for the proof, but outline the main calculation of the discrepancy between a weak
and a slic solution. Recall that
\begin{align*}
y(x,t) = w(R,t) \frac{x}{R} = t \varphi \Big ( \frac{R}{t} \Big ) \frac{x}{R}
\\
v ( R, t) = w_R  \left ( \frac{w}{R} \right )^{2}
\end{align*}
where $\varphi(s)$ is the self-similar solution of \eqref{mainsseq} in section \ref{sec-cavi}, and that
the specific volume  $v ( \cdot , t)$ is strictly monotone increasing
and obeys the bounds
$$
v_0 \le v(R,t) \le  \lambda^3 \, .
$$
The approximate solution  and the corresponding approximation of the specific volume are
\begin{align*}
y^n (x,t) = w^n (R,t) \frac{x}{R}
\\
w^n (R,t) =  \phi_n  \star w( \cdot , t)
\\
v^n = w^n_R \left ( \frac{w^n}{R} \right )^{2}
\end{align*}
where $\phi_n$ symmetric mollifier, $\phi (0) > 0$ and $w(\cdot , t)$ is an odd extension of $w$.
It can be proved that $v^n$ satisfies the bounds
\begin{align*}
c_1  &\le  v^n (R,t)  \le c_2  \qquad  \qquad \qquad \mbox{ for $R > \frac{1}{n}$}
\\
c_{\phi} n^3 w(0,t)^3 &\le v^n (R,t) \le c_3 (1 + t^3 n^3 )  \qquad \mbox{ for $R < \frac{1}{n}$  }
\end{align*}
where $c_{\phi}$ a positive constant depending on the shape of the mollifier.
Note that the approximate specific volume $v^n$ detects that there is a cavity forming, while
$v$ stays bounded away from zero.

To show we have a slic-solution we need to show that
$$
\frac{\del^2 y^n }{\del t^2} - \div \frac{\del W}{\del F} (\nabla y^n ) =: f^n \to 0 \quad \mbox{ in $\cD'$. }
$$
The discrepancy $D$  between the weak and the slic-solution lies in the
behavior near the cavity,  in the ball  $|x| < \frac{1}{n}$,
and  can be computed by the following heuristic "calculation"
\begin{align*}
D &= \iint S ( \nabla y^n) : \nabla \psi \;  \charf_{|x| < \frac{1}{n} } dx dt
\\
&\sim \int  \int_{|x| < \frac{1}{n} } \Big (  \frac{ w^n}{R} \Big )^{2} \, h' (v^n )  | \nabla \psi| \; R^{2} dR dt
\\
&\sim \int \int_0^1 \frac{ h' (n^3) }{n} \; | \nabla \psi| \;  \rho^{2} d\rho dt \, .
\end{align*}
Therefore, $D$ is related to the limit $\lim_{u \to \infty} \frac{h'(u^3) }{u}$. For the details and the
calculation of the discrepancy in the energy, see \cite{GT13}.

In conclusion, two issues may arise when a weak solution is viewed as limit of continuous solutions in a context
of strong singularities (like a cavity or fracture or shear band)

\begin{itemize}
\item[a.] New terms might appear in the momentum balance equation.

\item[b.] Even when such forces vanish in the limit, their effect can be felt in the energy balance and their contribution
might affect the total energy as the singularity forms.
\end{itemize}

\medskip
\noindent
{\bf Acknowledgements}
This research was supported by the EU FP7-REGPOT project "Archimedes Center for
Modeling, Analysis and Computation". AET is partially supported by the "Aristeia" program of the Greek
Secretariat for Research.


\begin{thebibliography}{10}


\bibitem{ball77}
{\sc J.M. Ball},
Convexity conditions and existence theorems in nonlinear elasticity,
{\em Arch. Rational Mech. Anal.} {\bf 63} (1977), 337-403.

\bibitem{BCO81}
{\sc J.M. Ball, J.C. Currie and P.J. Olver}
Null Lagrangians, weak continuity, and variational problems of arbitrary order
{\it J. Functional Analysis}  {\bf 41} (1981), 135-174.

\bibitem{ball82}
{\sc J.M. Ball},
Discontinuous equilibrium solutions and cavitation in
              nonlinear elasticity,
{\em Philos. Trans. Roy. Soc. London Ser. A}, {\bf 306}, (1982)  557--611.

\bibitem{Dafermos86}
{\sc C. Dafermos},
Quasilinear hyperbolic systems with involutions,
{\it Arch. Rational Mech. Anal.} {\bf 94} (1986), 373-389.




\bibitem{DST01}
{\sc S. Demoulini, D.M.A. Stuart, A.E. Tzavaras},
A variational approximation scheme for
three-dimensional elastodynamics with polyconvex energy,
{\em Arch. Rational Mech. Anal.} {\bf 157} (2001), 325-344.



\bibitem{edelen62}
{\sc D.G.B. Edelen},
The null set of the Euler-Lagrange operator
{\it Arch. Rational Mech. Anal.}  {\bf 11} (1962), 117-121.

\bibitem{ericksen62}
{\sc J.L. Ericksen},
Nilpotent energies in liquid crystal theories,
{\em Arch. Rational Mech. Anal.} {\bf 10} (1962), 189-196.

\bibitem{GT13}
{\sc J. Giesselmann and A.E. Tzavaras},
Singular limiting induced from continuum solutions and the problem of dynamic cavitation.
(submitted), (2013),
arXiv preprint arxiv:1306.6084.


\bibitem{MT12}
{\sc A. Miroshnikov and A.E. Tzavaras},
A variational approximation scheme for polyconvex elastodynamics that preserves the
positivity of Jacobians.
{\em Comm. Math. Sciences} {\bf 10} (2012), 87-115.


\bibitem{MT13}
{\sc A. Miroshnikov and A.E. Tzavaras},
On the construction and properties of weak solutions describing dynamic cavitation.
(preprint).

\bibitem{ps88}
{\sc K.A. Pericak-Spector and S.J. Spector},
Nonuniqueness for a hyperbolic system: cavitation in nonlinear elastodynamics.
{\em Arch. Rational Mech. Anal.} {\bf 101} (1988), 293 - 317.

\bibitem{ps97}
{\sc K.A. Pericak-Spector and S.J. Spector},
Dynamic cavitation with shocks in nonlinear elasticity.
{\em Proc. Royal Soc. Edinburgh Sect A} {\bf 127} (1997), 837 - 857.

\bibitem{qin}
{\sc T. Qin},
Symmetrizing nonlinear elastodynamic system,
{\em J. Elasticity} {\bf 50} (1998), 245-252.


\bibitem{SS03}
{\sc J. Sivaloganathan and S.J. Spector},
Myriad radial cavitating equilibria in nonlinear elasticity.
{\em SIAM J. Appl. Math.} {\bf 63} (2003), 1461 - 1473.

\bibitem{TN65}
{\sc C. Truesdell, W. Noll},
{\it The non-linear field theories of mechanics},
Handbuch der Physik {\bfseries III}, 3 (Ed. S.Fl\"{u}gge), Springer Verlag, Berlin, 1965.

\bibitem{wagner09}
{\sc D.H. Wagner},
Symmetric hyperbolic equations of motion for a hyper-elastic material,
{\it J. Hyper. Differential Equations} {\bf 6} (2009), 615-630.



\end{thebibliography}
\end{document}